\documentclass[11pt,a4paper]{amsart}
\usepackage{amssymb,amsmath}
\usepackage{latexsym}
\usepackage{xypic}

\newcommand{\F}{\mathbb F}
\newcommand{\PP}{\mathbb P}

\newcommand{\kO}{\mathcal{O}}

\newcommand{\cV}{\mathcal{V}}
\newcommand{\cP}{\mathcal{P}}

\newtheorem{theorem}{Theorem}[section]

\newtheorem{coro}[theorem]{Corollary}

\newtheorem{rem}[theorem]{Remark}

\newcommand{\prf}{\noindent{\em Proof}. }

\newcommand{\T}{{\mathbb T}}
\newcommand{\C}{{\mathbb C}}
\newcommand{\N}{{\mathbb N}}

\newcommand{\Q}{{\mathbb Q}}
\newcommand{\R}{{\mathbb R}}
\newcommand{\Z}{{\mathbb Z}}

\newcommand{\ra}{\rangle}
\newcommand{\la}{\langle}

\newcommand{\Pic}{{\rm Pic}}
\newcommand{\Hom}{{\rm Hom}}

\newcommand{\verylongarrow}[1]{\hbox to #1{\rightarrowfill}}

\newcommand\mynote[1]{\marginpar{\ \\ \small \tt #1}}
\newcommand\bel[1]{{\mynote{#1}}\begin{equation}\label{#1}}

\begin{document}
\title[Integral Cohomology and Mirror Symmetry]{
Integral Cohomology and Mirror Symmetry for
\\Calabi-Yau 3-folds}

\author[Victor Batyrev]{Victor Batyrev}
\address{Mathematisches Institut, Universit\"at T\"ubingen, Auf der 
Morgenstelle 10, 72076 T\"ubingen, Germany}
\email{victor.batyrev@uni-tuebingen.de}

\author[Maximilian Kreuzer]{Maximilian Kreuzer}
\address{Institut f\"ur Theoretische Physik,
Technische Universit\"at Wien,
Wiedner Hauptstr. 8-10/136,
A-1040 Vienna, Austria} 
\email{maximilian.kreuzer@tuwien.ac.at}

\begin{abstract}
In this paper we compute the integral cohomology groups 
for all examples of Calabi-Yau 
3-folds obtained from  
hypersurfaces in $4$-dimensional Gorenstein toric Fano varieties. 
Among 473 800 776 families of Calabi-Yau 3-folds $X$  
corresponding to $4$-dimensional reflexive polytopes there exist 
exactly 32 families having non-trivial torsion  
in $H^*(X, \Z)$. We came to an interesting observation that 
the torsion subgroups in $H^2$ and $H^3$ are exchanged by 
the mirror symmetry involution. 
\end{abstract}

\maketitle
\section*{Introduction}

Let $X$ be  a smooth projective Calabi-Yau $d$-fold over $\C$.  
We shall always assume that $h^i(X, {\kO}_X) =0$ for $0 < i <d$. 
Our main interest 
is the torsion in the integral cohomology ring $H^*(X, \Z)$.
For a finite abelian group $G$ we denote by $G^*$ the 
dual abelian group $\Hom (G, \Q/\Z)$.
 Using the isomorphism from the universal coefficient theorem
\[  {\rm Tors}(H_i(X, \Z))\cong  [ {\rm Tors}(H^{i+1}(X, \Z))]^* \]
and the  Poincar\'e duality 
\[  {\rm Tors}(H_i(X, \Z))\cong  {\rm Tors}(H^{2d-i}(X, \Z)),\]
we obtain for $d=3$:   
\[ H^0(X, \Z) \cong   H^6(X, \Z) \cong \Z,\;\; H^1(X, \Z) = 0, \]
\[ H^2(X, \Z) \cong A(X) \oplus \Z^a, \;\;  
H^4(X, \Z) \cong  B(X)^* \oplus \Z^a \]
\[ H^3(X, \Z) \cong  B(X) \oplus \Z^{2b +2}, \;\; H^5(X, \Z) \cong  A(X)^*,\]
where the finite abelian groups
\[ A(X):= {\rm Tors}(H^2(X, \Z)), \;  B(X):= {\rm Tors}(H^3(X, \Z)) \]
determine completely the  torsion in $H^*(X, \Z)$.
In particular, we have the  isomorphisms 
\[ A(X) \oplus B(X)^* \cong {\rm Tors}( H^{\rm even}(X, \Z))  \cong  
[ {\rm Tors}(H^{\rm odd}(X, \Z))]^*. \]
It is known that mirror symmetry exchanges the integers  
$a,b$ so that  for a mirror
Calabi-Yau $3$-fold $X^*$ one has
\[ \Q^{2a+2} \cong H^{\rm even}(X, \Q)  \cong  H^{\rm odd}(X^*, \Q), \]
\[ \Q^{2b+2} \cong H^{\rm odd}(X, \Q)  \cong  H^{\rm even}(X^*, \Q). \]
An interesting mathematical question is to understand the 
behavior of the torsion in $H^*(X, \Z)$, i.e. the groups $A(X)$ and
$B(X)$, under the mirror symmetry.  
 The group
$B(X)$ is isomorphic to the torsion in the \'etale cohomology group 
$H^2_{\mathrm{et}}(X, {\kO}^*_X)$ 
and it  is called the {\em cohomological  
Brauer group of} $X$. By a recent result 
of Kresch and Vistoli \cite{KV}, $B(X)$  is also isomorphic to the 
Brauer group of Azumaya algebras on $X$. On the other hand, there are  canonical isomorphisms 
$$A(X) \cong \Hom (\pi_1(X), \Q/\Z) \cong {\rm Tors} ({\rm Pic}(X)).$$

From the physical point of view 
it is more natural to consider the topological $K$-groups  
$K^0(X)$ and $K^1(X)$ together with  
some natural filtrations \cite{BD,BDM}. 
Using the Atiyah-Hirzebruch spectral sequence (see \cite{BD}, 2.5), 
one obtains two short exact sequences
\[ 0 \to A(X)^* \to  {\rm Tors}(K^1(X)) \to B(X) \to 0, \]
\[  0 \to B(X)^* \to  {\rm Tors}(K^0(X)) \to A(X) \to 0. \]
It is expected that the mirror symmetry exchange $K^0$ and $K^1$ so that for a mirror 
pair $(X, X^*)$ one has 
the isomorphisms 
\[  {\rm Tors}(K^1(X)) \cong {\rm Tors}(K^0(X^*)), \;\; 
  {\rm Tors}(K^0(X)) \cong {\rm Tors}(K^1(X^*)). \] 
 These isomorphisms agree with predictions of SYZ-construction 
and topological calculations of M. Gross \cite{Gross}. 
The compatibility of the above isomorphisms 
with the natural filtrations in $K^i$ would imply 
the isomorphisms 
\begin{equation} A(X) \cong  B(X^*) , \;\; B(X) \cong A(X^*). 
\label{tor-dual}
\end{equation}
The main purpose of this paper is to verify these isomorpisms for 
all examples of Calabi-Yau 3-folds obtained from  
hypersurfaces in $4$-dimensional Gorenstein toric Fano varieties.

There exist 473 800 776 families of Calabi-Yau 3-folds $X$  
corresponding to $4$-dimensional reflexive polytopes \cite{KS1,KS2}. 
These families give rise to 30 108 different pairs of numbers 
$(a,b)$.
We show that among 473 800 776 families of Calabi-Yau 3-folds $X$  
there exist 
exactly 32 families having non-trivial torsion in $H^*(X, \Z)$.
More precisely, there are exactly 16 families of simply-connected
Calabi-Yau $3$-folds $X$ having 
non-trivial Brauer group $B(X) \cong 
\Z/p\Z$ $(p =2,3,5)$. They are mirror dual 
to another  16 families of Calabi-Yau $3$-folds $X^*$ having  trivial 
Brauer group and a non-trivial cyclic fundamental group of order 
$p =2,3,5$ \cite{KKRS}. Thus, we come to the observation that for all 
families of Calabi-Yau 
3-folds  one has the isomorphisms (\ref{tor-dual}). Although 
the groups $A(X)$ and $B(X)$ can be computed in a purely combinatorial way 
using lattice points in faces of $4$-dimensional 
reflexive polytopes $\Delta$ 
(or dual reflexive polytopes $\Delta^*$) we do not see immediately 
the isomorphisms  (\ref{tor-dual}) from the combinatorial duality 
between $\Delta$ and $\Delta^*$. Therefore, 
it would be interesting to find a general mathematical 
explanation for the mirror symmetry isomorphism
\[ B(X) \cong \Hom (\pi_1(X^*), \Q/\Z). \]

If  a finite group $G$ acts on a smooth Calabi-Yau manifold $V$, then  
there exists  a conformal field theory of   the orbifold $V/G$  associated 
with an element $\alpha \in H^2(G, U(1))$ (discrete torsion).  
In \cite{VW} Vafa and Witten have suggested  
that there should be some connection between elements in the Brauer group 
of $V/G$ ($B$-fields) and elements  in $H^2(G, U(1))$ 
(discrete torsion
$2$-cocycles).  
In this paper we show that all 16 families of Calabi-Yau $3$-folds 
with non-trivial cyclic 
Brauer group are birational to Calabi-Yau orbifolds $V/G$ 
($G \cong \Z/p\Z \oplus \Z/p\Z$, $ p =2,3,5$) 
and the corresponding elements in the cyclic 
Brauer group $B(X) \cong H_2(G) \cong\Z/p\Z$ 
can be identified with the discrete torsion  
$2$-cocycles.  

It is easy to see that 
there are  no algebraic curves in $X$ 
representing the torsion classes in $H_2(X, \Z)$.
Aspinwall and Morrison  suggested in  \cite{AM2} that the instanton 
expansion of the Yukawa coupling for a Calabi-Yau 3-fold with 
the Brauer group $B(X) \cong \Z/p\Z$ 
should  be weighted by powers of $p$-th roots of 
unity. It would be interesting to verify this prediction for the 
above 16 toric families 
of Calabi-Yau 3-folds with non-trivial Brauer group.
\medskip

\noindent
{\bf Acknowledgements:} We would like to thank the  
Fields Institute for 
Research in Mathematical Sciences in Toronto 
for hospitality during the spring of 2005 and
DFG-Schwerpunkt ``Global methods in complex geometry'' for the support
of the first author. The second author was supported by FWF 
grant Nr. P15584-N08. 
We are grateful to Noriko Yui for her 
interest and encouragement in our work. We thank also Volker Braun, 
Mark Gross, Kentaro Hori, Klaus Hulek, Askold Khovanski, Tony Pantev,  
Duco van Straten, Bernd Siebert and  
Bal{\'a}zs\ Szendr\H{o}i for useful discussions. 
\bigskip

\section{The fundamental group of toric 
hypersurfaces}\label{fund}

First we recall a combinatorial computation of the fundamental group 
of a toric variety \cite{D}. Let $M \cong \Z^d$ be a lattice of rank $d$ and 
$N:= \Hom (M, \Z)$ the dual lattice. We denote by 
$\la *, * \ra$ the canonical pairing $M \times N \to \Z$. 
For an algebraic torus 
$\T^d := {\rm Spec}\, \C[M] \cong (\C^*)^d $ there are canonical isomorphisms
\[ \pi_1^{\rm top}(\T^d) \cong N \cong  H_1(\T^d, \Z), \;\; 
M \cong  H^1(\T^d, \Z). \]
If $\PP_{\Sigma}$ is a smooth partial toric compactification of $\T^d$ defined 
by a fan $\Sigma \subset N_{\R} = N \otimes \R$, then 
$\PP_{\Sigma} \setminus \T^d$ is a union of divisors $D_1, \ldots, D_n$ which 
1-to-1 correspond to primitive lattice vectors $e_1, \ldots, e_n \in N$ generating 
$1$-dimensional cones in $\Sigma$.  
If $D_i^{\circ}$ is a dense torus orbit in $D_i$, then 
we have $\PP_i:= \T^d \cup D_i^{\circ} \cong \C \times \T^{d-1}$. 
The embedding $\T^d \hookrightarrow \PP_i$ defines 
the homomorphism of the fundamental groups 
\[ \rho_i\, : \, N = \pi_1^{\rm top} (\T^d) \to 
\pi_1^{\rm top}(\PP_i) \cong  N/\Z e_i.  \]
Combining 
all the 
homomorphisms $\rho_i$ $(i =1, \ldots, n)$, we obtain 
the isomorphisms  
\[  \pi_1^{\rm top}(\PP_{\Sigma})  \cong  
\pi_1^{\rm top}(\PP_{\Sigma^{(1)}}) \cong  N/\sum_{i=1}^n \Z e_i,  \]
where $\PP_{\Sigma^{(1)}}$ is a toric variety associated to 
the subfan $\Sigma^{(1)} \subset \Sigma$ of all cones of dimension $\leq 1$ in 
$\Sigma$. In particular, we have 
\[ H_1(\PP_{\Sigma}, \Z) \cong   H_1(\PP_{\Sigma^{(1)}}, \Z) \cong 
 N/\sum_{i=1}^n \Z e_i.  \]

In oder to compute the fundamental group of hypersurfaces in toric varieties
one applies the following result of Oka \cite{Oka}:

\begin{theorem} Let $\Delta \subset M_{\R}$ be a lattice polytope of dimension 
$d$,  
$$f(t) = \sum_{m \in \Delta\cap M} c_m t^m \in 
\C[M] \cong \C[ t_1^{\pm 1}, \ldots, 
t_d^{\pm 1} ]$$
be a $\Delta$-nondegenerate 
Laurent polynomial, and $W \subset \T^d$ the nondegenerate 
affine hypersurface defined by 
the equation $f(t) =0$. Then the embedding $j\, : \,  W 
\hookrightarrow \T^d$ is 
a homotopy $(d-1)$-equivalence, i.e. the 
homomorphisms of the homotopy groups 
\[  j_* \; : \; \pi_i(W) \to \pi_i(\T^d) \]
is bijective for $i < d-1$ and surjective for $i = d-1$. 
\end{theorem} 

By Whitehead's theorem, we get: 

\noindent
\begin{coro} The homomorphism 
\[   j_* \; : \; H_i(W, \Z) \to H_i(\T^d, \Z) \]
is bijective for $i < d-1$ and surjective for $i =d-1$. 
\end{coro} 

Since $H_i(\T^d, \Z)$ is a free abelian group of rank ${ d \choose i }$, 
using the universal coefficients theorem,  we also obtain: 

\noindent
\begin{coro} For any abelian group ${\mathcal A}$,  the homomorphism 
\[   j^* \; : \;  H^i(\T^d, {\mathcal A}) \cong  \Lambda^i M \otimes_{\Z} 
{\mathcal A} 
\to H^i(W, {\mathcal A})   \]
is  bijective for $i < d-1$ and injective for $i =d-1$. 
\label{L-coh}
\end{coro} 

Let $\Sigma:= \Sigma(\Delta) \subset N_{\R}$ be the normal fan 
of the polytope $\Delta \subset M_{\R}$. The cones $\sigma= \sigma({\theta}) 
\in 
\Sigma({\Delta})$ 1-to-1 correspond to faces $\theta \subset \Delta$ in the 
following way:
\[ \sigma({\theta}) = \{ y \in N_{\R}\; : \; \min_{ x \in \Delta}
\la x, y \ra = \la x, z \ra \;\; \forall z \in \theta \}, \]
i.e. $\sigma({\theta})$ consists of those $y \in N_{\R}$ for which 
the minimum of the linear function $\la *, y \ra$ on $\Delta$ 
is attained at any point of  $\theta$. Sometimes we denote the corresponding 
projective toric variety $\PP_{\Sigma} = 
\PP_{\Sigma({\Delta})}$ simply by $\PP_{\Delta}$.
Let $\overline{W} \subset \PP_{\Sigma}$ be the Zariski closure 
of $W$ in   $\PP_{\Sigma}$. If $\Sigma^{(1)} \subset \Sigma$ 
is the subfan of all cones of dimension $\leq 1$ in $\Sigma$, then 
 $\PP_{\Sigma^{(1)}}$ and  $\overline{W}^{(1)} := \PP_{\Sigma^{(1)}} \cap 
\overline{W}_f$ are  smooth quasi-projective varieties. 

\noindent
\begin{theorem} 
Let $d \geq 3$. Then the embedding $j \, : 
\,\overline{W}^{(1)} \hookrightarrow  
 \PP_{\Sigma^{(1)}}$ induces  the isomorphism
\[  j_*\; : \; \pi_1(\overline{W}^{(1)}) \cong \pi_1( \PP_{\Sigma^{(1)}}) \cong 
N/\sum_{i =1}^n \Z e_i,  \]
where $e_1, \ldots, e_n$ are primitive generators of $1$-dimensional 
cones in $\Sigma^{(1)}$.  
\label{w1}
\end{theorem}

\noindent
\prf Let $ \PP_{\Sigma^{(1)}} \setminus \T^d = D_1^{\circ} \cup 
\cdots  D_n^{\circ}$ $(  D_i^{\circ} \cong \T^{d-1}$). 
Using the isomorphisms  $\pi_1(W) \cong \pi_1( \T^d) \cong N$ and 
the surjectivity of $\pi_1(W) \to  \pi_1(\overline{W}^{(1)})$, we obtain 
that 
\[ j_*\, : \, \pi_1(\overline{W}^{(1)}) \to \pi_1( \PP_{\Sigma^{(1)}}) \]
is surjective. Therefore the universal covering $\overline{U}^{(1)}$ of 
$\overline{W}^{(1)}$ is induced by some unramified covering $U$ of 
$\T^d$.

In order to prove that $j_*$ is an isomorphism it remains
to show that the image of every element $e_i \in\pi_1(W) \cong N$ 
($i =1, \ldots, n)$ is zero in  $\pi_1(\overline{W}^{(1)})$. Let  
$\PP_i:=\T^d \cup D_i^{\circ}$ and $\overline{W}_i^{(1)}: =  
\overline{W}^{(1)} \cap \PP_i$.
Using open inclusions 
$W \subset \overline{W}^{(1)}_i \subset 
 \overline{W}^{(1)}$, it is enough to prove that 
the image of $e_i \in\pi_1(W)$ is zero already 
in $\pi_1( \overline{W}^{(1)}_i)$
so that one obtains  $\pi_1( \overline{W}^{(1)}_i) \cong N/\Z e_i$. We
can consider $ \overline{W}^{(1)}_i$ as an affine hypersurface in 
$\PP_i \cong \T^{d-1} \times \C \cong  {\rm Spec}\, \C[t_1^{\pm 1}, 
\ldots, t_{d-1}^{\pm 1}, t_d]$ where the divisor 
$D_i^{\circ} \subset \PP_i$ is defined by $t_d =0$.
Moreover, the element $e_i \in \pi_1(\T^d) = N$ is represented by 
a small 1-cycle $\gamma_i \subset \C^*$ around $0 \in \C$ 
which can be contracted to $0$ in $\T^{d-1} \times \C$.    
The transversality of $ \overline{W}^{(1)}_i$ and $D_i^{\circ}$ 
implies that the 1-cycle $\gamma_i$ can be choosen inside 
$W_f$ so that it will be contractible in  $ \overline{W}^{(1)}_i$. 
Therefore all elements $e_1, \ldots , e_n$ are in the kernel 
of the homomorphism  $\pi_1(W) \to  \pi_1(\overline{W}^{(1)})$.
We remark that the isomorphism $\pi_1( \overline{W}^{(1)}_i) \cong 
\pi_1(\PP_i)$ can be interpreted also as a Lefschetz-type statement 
for a nondegenerate hypersurface $\overline{W}^{(1)}_i$ 
in $\T^{d-1} \times \C$.
\hfill $\Box$ 
\medskip

Since $e_1, \ldots, e_n$ generate $N$ over $\Q$, we obtain:

\noindent
\begin{coro} 
The fundamental group of the quasi-projective hypersurface 
$\overline{W}^{(1)}$ is a finite abelian group. In particular,
$ \pi_1(\overline{W}^{(1)})$ coincides with the algebraic
fundamental group of $\overline{W}^{(1)}$ and the universal 
covering of  $\overline{W}^{(1)}$ is a toric quasi-projective 
hypersurface $\widetilde{W}^{(1)}$ 
in the finite universal covering of $\PP_{\Sigma^{(1)}}' \to
 \PP_{\Sigma^{(1)}}$ 
determined by the finite index sublattice 
$N' = \sum_{i=1}^n \Z e_i \subset N$.
\label{prep-fund}
\end{coro} 

\begin{theorem} 
Assume that $d \geq 3$. 
Let $\PP_{\widehat{\Sigma}}$ be a smooth $d$-dimensional 
projective toric variety 
defined by a fan $\widehat{\Sigma} \subset N_{\R}$ which is a 
subdivision of $\Sigma= \Sigma(\Delta)$. We denote by $\Gamma$ the 
set of primitive generators $v$ 
of $1$-dimensional cones in 
 $\widehat{\Sigma}^{(1)}$ such that the 
minimum of $\la *, v \ra$ on $\Delta$ is attained 
on a face $\theta \subset \Delta$ of dimension $\geq 1$. 
Define the sublattice $N_{\Delta}^{(1)} 
\subset N$ to be generated by all $v \in \Gamma$. 
If $\widehat{W}$ is the closure of $W$
in  $\PP_{\widehat{\Sigma}}$, then the fundamental group $\pi_1(\widehat{W})$ 
is isomorphic to the cyclic group  
\[ N/ N_{\Delta}^{(1)}. \]
\label{fundam}
\end{theorem} 

\prf  Let ${\mathcal E} := 
\{e_1, \ldots, e_n \}$ be the set of primitive generators 
of $1$-dimensional cones in 
${\Sigma}^{(1)}$. Since every linear function $\la *, e_i \ra$ attains 
its minimum of $\Delta$ on a $(d-1)$-dimensional face of $\Delta$, 
we have  ${\mathcal E} \subset \Gamma$ and we can write 
$\Gamma = \{ e_1, \ldots, e_n, \ldots, e_{n+k} \}$. Moreover, 
 $\overline{W}^{(1)}$ can be considered as a dense 
open subset of  $\widehat{W}$. Therefore, we have a surjective  
homomorphism 
\[ \psi\; : \; \pi_1(\overline{W}^{(1)}) \to \pi_1(\widehat{W}) \]
and the dual injective homomorphism 
\[ \psi^* \; : \; \Hom (  \pi_1(\widehat{W}, \Q/Z) \to 
  \Hom (  \pi_1(\overline{W}^{(1)}),  \Q/Z). \]
Since $\widehat{W} \setminus \overline{W}^{(1)} = Z_1 \cup \cdots Z_k$, where 
$$Z_i = D_{n+i} \cap \widehat{W}, \;\; i =1, \ldots, k,$$
we obtain that an element $g \in  \Hom (  \pi_1(\overline{W}^{(1)}),  \Q/Z) 
\cong \Hom ( N/\sum_{i=1}^n \Z e_i, \Q/\Z)$  
belongs to ${\rm Im}\, \psi^*$ if and only if the cyclic Galois 
covering of $\overline{W}^{(1)}$ corresponding to $g$ is unramified
along  all divisors $Z_1, \cdots, Z_k$. Using the same arguments as 
in the proof of \ref{w1}, we obtain that the latter holds exactly 
when $g(e_{n+1}) =  \cdots =g(e_{n+k}) =0$. Therefore, 
\[ {\rm Im}\, \psi \cong  \Hom (  N/\sum_{i=1}^{n+k} \Z e_i, \Q/\Z) = 
 \Hom (  N/ N_{\Delta}^{(1)}, \Q/\Z) \]
and by duality $  \pi_1(\widehat{W}) \cong  N/ N_{\Delta}^{(1)}$. 
The group  $N/ N_{\Delta}^{(1)}$ is cyclic, because the smoothness of 
$\widehat{W}$ implies that for any $(d-1)$-dimensional cone 
$\sigma = \sigma(\theta) \in \Sigma$ the subset $\sigma \cap \Gamma$ 
generate a direct summand of the rank $d-1$ of the lattice $N$ 
(all $(d-1)$-dimensional cones in $\widehat{\Sigma}$ which subdivide
$\sigma$ are generated by part of a $\Z$-basis of $N$). 
\hfill $\Box$

\begin{coro} 
For any lattice point $v \in N$, we denote by 
$$M_v := \{ x \in M \; : \; 
\langle x, v \rangle = 0 \}$$ 
the orthogonal complement to $v$ in 
$M$. Then we have the exact sequence
\[   
\bigoplus_{v \in \Gamma} \Lambda^{d-1} M_v \to 
  \Lambda^{d-1} M \to H^{2d-3}_{\rm tor}(\widehat{W},\Z) \to 0. \]
\end{coro}

\prf The exactness of 
 \[  
\bigoplus_{v \in \Gamma} \Lambda^{d-1} M_v \to 
  \Lambda^{d-1} M \to H^{2d-3}_{\rm tor}(\widehat{W},\Z) \to 0 \]
follows from the exactness of 
\[ \bigoplus_{v \in \Gamma} \Z [ v] \to N \to  \pi_1(\widehat{W}) \to 0 \]
together with the canonical isomorphisms
\[  \pi_1(\widehat{W})  \cong  H^{2d-3}_{\rm tor}(\widehat{W},\Z), 
\;\;  \Lambda^{d-1} M_v \cong \Z, \;\; \Lambda^{d-1} M \cong N \]
after choosing 
an isomorphism (orientation) 
$\Lambda^d \cong \Z$.
\hfill $\Box$

\begin{rem} 
{\rm Theorem \ref{fundam} can be considered as a special case of a 
more general result of Oka (see \cite{Oka-b}, Chapter V, \S 5). } 
\end{rem}  

\noindent
\begin{coro} 
Let $\overline{W}$ be a nondegenerate 
$(d-1)$-dimensional Calabi-Yau hypersurface in the Gorenstein 
toric Fano variety $\PP_{\Delta}$ associated with a $d$-dimensional 
reflexive polytope $\Delta \subset M_{\R}$. 
Assume 
$d \geq 3$ and 	that
there exists a smooth projective crepant 
resolution $\widehat{W}$ of singularities in  $\overline{W}$ induced by  
a convex (coherent) triangulation of the dual reflexive polytope $\Delta^* 
\subset \N_{\R}$. We denote by  
$N_{\Delta^*}'$ the sublattice in $N$ generated by all lattice 
points in $\Delta^* \cap N$ which belong to faces of codimension 
$>1$ of $\Delta^*$. Then the fundamental group of  $\widehat{W}$
is isomorphic to the cyclic group 
$$N/N_{\Delta^*}'.$$
\label{toric-h-p1}
\end{coro}

\prf The statement follows immediately from \ref{fundam}, because for 
a crepant desingularization $\widehat{W} \to \overline{W}$ 
the set of primitive 
generators $\Gamma$ is exactly the set of all lattice 
points in faces of codimension $<1$ of 
the dual reflexive polytope $\Delta^*$ (see \cite{B1}). 
\hfill $\Box$

\begin{rem} 
{\rm Recently Haase and  Nill  proved that for a reflexive polytope
$\Delta^*$ of arbitrary dimension $d > 2$ 
the sublattice in $N$ generated by {\bf all} lattice points
in $\Delta^* \cap N$ coincides with $N_{\Delta^*}'$ 
\cite{HN}. The proof of this results 
uses the fact that the interior lattice points in codimension-$1$ 
faces of $\Delta^*$ are Demazure roots for 
the authomorphism group of the corresponding Gorenstein 
toric Fano variety $\PP_{\Delta^*}$ \cite{Nill1}. In particular, 
this result easily implies that $N = N_{\Delta^*}'$ 
for all reflexive polyhedra $\Delta^*$ of dimension 
$d =3$. The last statement follows also from 
\ref{toric-h-p1} and from the fact 
that every smooth $K3$-surface is simply connected.}
\end{rem} 

\begin{rem} 
{\rm It is known that for $d =4$ there always 
exists a smooth projective crepant 
resolution $\widehat{W}$ of singularities in  $\overline{W}$ induced by  
a convex (coherent) triangulation of the dual reflexive polytope $\Delta^* 
\subset \N_{\R}$. For $d \geq 5$ such a resolution does not exist in general.
However, one can use $N/N_{\Delta^*}'$ as a candidate for a
{\bf stringy fundamental group} of a singular Calabi-Yau variety 
$\overline{W}$ or its maximal partial crepant desingularization  
$\widehat{W}$ (similar to stringy Hodge numbers 
\cite{B2}). It would be interesting  to know whether  for $d \geq 5$ 
the dual group  $(N/N_{\Delta^*}')^*$ coincides 
with the torsion subgroup in the Picard group 
of the Deligne-Mumford Calabi-Yau stack associated with the $V$-manifold
$\widehat{W}$ (see \cite{BCS,PS1,PS2}).  
}
\end{rem}

\section{The list of 
non-simply connected Calabi-Yau 3-folds}\label{ex}
  
Using the Calabi-Yau data \cite{KS2} and the program package PALP \cite{PALP} 
one can check all 473 800 776 reflexive
polytopes $\Delta$ and find that there exist exactly 16 examples
of reflexive polytopes 
$\Delta_i$ $(1 \leq i \leq 16)$ such that $|N/N_{\Delta^*_i}'|>1$. 
In the latter case, 
 $|N/N_{\Delta^*_i}'|$ is always a prime number $p = 2,3,5$.  

The most known  example of a  non-simply connected Calabi-Yau 3-fold 
obtained from a hypersurface in a Gorenstein toric variety $\PP_{\Delta}$
is a free $\mu_5$-quotient of a smooth quintic 3-fold in $\PP^4$ 
obtained as follows. 
 Consider the  action of  the cyclic group $\mu_5 = \la \zeta \ra$
 by $(1, \zeta, \zeta^2,\zeta^3,\zeta^4)$ on $\PP^4$ and take  
$\PP_{\Delta_1} := \PP^4/\mu_5$. The Gorenstein toric Fano 
variety $\PP_{\Delta_1}$ contains exactly 5 isolated singular points.
A non-degenerate hypersurface $\overline{W}_1 \subset \PP_{\Delta_1}$
is smooth, i.e. $\overline{W}_1 = \widehat{W}_1$, 
it has the Hodge numbers $h^{1,1}(\overline{W}_1) =
 1$, $h^{2,1}(\overline{W}_1) = 21$, and $\pi_1(\overline{W}_1) \cong \Z/5\Z$.
This is the single example in the case $|N/N_{\Delta^*_i}'|=5$.

Similarly one obtains 
 a free $\mu_3$-quotient $\overline{W}_2$ 
 of a smooth hypersurface of bidegree $(3,3)$ 
in $\PP^2 \times \PP^2$ using the action of  the cyclic 
group $\mu_3 = \la \zeta \ra$ on $\PP^2 \times \PP^2$ by 
$(1, \zeta, \zeta^2)$ on each $\PP^2$. 
The toric variety $\PP_{\Delta_2} := \PP^2 \times 
\PP^2/\mu_3$ has $9$ isolated 
singularities. A non-degenerate hypersurface 
$\overline{W}_2 = \widehat{W}_2  \subset \PP_{\Delta_2}$
is smooth, it  has the Hodge numbers $h^{1,1}(\overline{W}_2) =
 2$, $h^{2,1}(\overline{W}_2) = 29$, and 
$\pi_1(\overline{W}_2) \cong \Z/3\Z$.

We remark 
that the dual reflexive polytopes $\Delta_1^*, \ldots, \Delta_{16}^*$ 
have at most $8$ vertices. $\Delta_{14}^*$ is the single dual reflexive 
polytope with $8$ vertices. The corresponding Gorenstein 
toric Fano variety $\PP_{\Delta_{14}}$ is  $\PP^1 \times \PP^1 
\times \PP^1 \times \PP^1/\mu_2$ 
where the cyclic group $\mu_2 = \la \zeta \ra$ acts by 
$(1, \zeta)$ on each $\PP^1$ . 
\bigskip 

We collected the information about all 16 reflexive polytopes 
$\Delta_1, \ldots, \Delta_{16}$ in the table below. 
In this  table  
$\cP_{\Delta}=|\Delta\cap M|$, $\cV_{\Delta}=|\mathrm{vert}(\Delta)|$,
        $\cP_{\Delta^*}=|\Delta^*\cap N|$, and $\cV_{\Delta^*}=
|\mathrm{vert}(\Delta^*)|$ 
denote numbers of points and vertices, respectively. 
Since any reflexive polytope $\Delta_i \subset M_{\R}$ is uniquely determined
by its dual polytope $\Delta^*_i \subset N_{\R}$, we prefer to use 
the small number  of vertices $v_1, \ldots, v_{\cV_{\Delta^*}}$ 
 in $\Delta_i^*$ 
in order to describe $\Delta_i^*$ by $\cV_{\Delta^*} -4$ 
independent linear relations 
among them. 

The lattice $N$ is generated by the vertices  
$v_1, \ldots, v_{\cV_{\Delta^*}}$ and by an additional vector 
$v$ expressed as a rational linear combination of the vertices $\{v_i \}$.
The sublattice $N_{\Delta_i^*}'$ is generated by 
$v_1, \ldots, v_{\cV_{\Delta^*}}$ except for 4 cases: 
$i = 4,9, 15,16$. In the latter cases  $N_{\Delta_i^*}'$ is generated 
by $v_1, \ldots, v_{\cV_{\Delta^*}}$ and by the vector $2v$. 

\newpage

\begin{center}

\begin{tabular}{||c|c|c@{~}c|c@{~}c|c@{~}c|c|c||}
                                                        \hline\hline
$n^{\circ}$&$|\pi_1|$ & $\cP_{\Delta}$ & $\cV_{\Delta}$ & $\cP_{\Delta^*}$ 
& $\cV_{\Delta^*}$ & 
        $h^{1,1}$ & $h^{2,1}$ & $\chi(\widehat{W})$ & 
        vertices  $\{v_i\}$  of  ${\Delta}^*$, $N:= 
\Z v + \sum_i \Z v_i$ 
                              \\\hline\hline
1&$5$ & 26 & 5& 6 & 5 & 1 & 21 & -40&    
$ \begin{array}{c} v_1 + v_2 + v_3 + v_4 + v_5 =0 \\
v = \frac15(v_2 + 2 v_3 + 3v_4 + 4 v_5) \end{array}$ 
                                                        \\\hline\hline

2&$3$& 34&9& 7&6& 2&29&-54& $ \begin{array}{c}
 v_1 + v_2 + v_3 = v_4 + v_5 + v_6 =0 \\  
 v = \frac13(v_2 + 2 v_3 + v_5 + 2v_6) \end{array} $ \\\hline
3& $3$& 49&5& 7&5& 2&38&-72&   $\begin{array}{c}
 3v_1 + 3v_2 + v_3 + v_4 + v_5 =0 \\
v = \frac13(v_1+ 2v_2 + v_3 + 2v_4) \end{array}$
                          \\\hline\hline
4&$2$& 53 & 5 & 9 &5 &3&43 &-80&  $\begin{array}{c}
 4v_1 + v_2 + v_3 + v_4 + v_5 =0 \\
v =  \frac14(2v_1+ v_2 + 2v_3 + 3v_4) \end{array}$ \\\hline
5&$2$& 77&7& 9&6& 3&59&-112&  $\begin{array}{c}
 4v_1 + 2v_2 + v_3 + v_4 =0 \\
 2v_1 + v_5 + v_6 =0 \\
v = \frac12(v_1+ v_2 + v_3 + v_5) \end{array}$
 \\\hline
6&$2$& 77&9& 9&7& 3&59&-112& $\begin{array}{c}
 2v_1 + v_2 + v_3 = 2v_1 + v_4 + v_5 =0 \\
 2v_1 + v_6 + v_7 =0 \\
v = \frac12(v_1+ v_2 + v_4 + v_6) \end{array}$ \\\hline
7&$2$& 101&5& 9&5& 3&75&-144&  $\begin{array}{c}
 8v_1 + 4v_2 + 2v_3 +v_4 + v_5 =0 \\  
 v =\frac12 (v_1 + v_2 +  v_3 + v_4) \end{array} $ 
            \\\hline
8&$2$& 101&6& 9&6& 3&75&-144& 
$ \begin{array}{c}
 4v_1 + 2v_2 + v_3 + v_4 =0 \\
 4v_1 + 2v_2 + v_5 + v_6 =0 \\  
 v = \frac12(v_1 + v_2 + v_3 + v_5) \end{array}$
                                \\\hline
9&$2$& 29&8& 9&6& 4&28&-48& $\begin{array}{c}
 v_1 + v_2 + v_3 + v_4 = v_5 + v_6 =0 \\
v =  \frac14(v_2+ 2v_3 + 3v_4 + 2v_6) \end{array}
$ 
                    \\\hline
10&$2$& 53&8& 9&6& 4&44&-80& 
 $\begin{array}{c}
 4v_1 + 2v_2 + v_3 + v_4 =v_5 + v_6 =0 \\
v = \frac12(v_1 + v_2+ v_3 + v_5) \end{array}
$                  
                        \\\hline
11&$2$& 53&10& 9&7& 4&44&-80&  $\begin{array}{c}
 2v_1 + v_2 + v_3 =0 \\
 2v_1 + v_4 + v_5 =0 \\
 v_6 + v_7 =0 \\
v =\frac12 (v_1+ v_2 + v_4 + v_6) \end{array}$
                        \\\hline
12&$2$& 41&9& 9&6& 4&36&-64&   $\begin{array}{c}
 2v_1 + v_2 + v_3 =0 \\
 2v_4 + v_5 + v_6 =0 \\
v = \frac12(v_1+ v_2 + v_4 + v_5) \end{array}$\\\hline
13&$2$& 41&12& 9&7& 4&36&-64& 
 $\begin{array}{c}
 2v_1 + v_2 + v_3 =0 \\
 v_4 + v_5 =  v_6 +v_7 =0 \\
v = \frac12(v_1+ v_2 + v_4 + v_6) \end{array}$ \\\hline
14&$2$& 41&16& 9&8& 4&36&-64&  $\begin{array}{c}
 v_1 + v_2 = v_3 +v_4 =0 \\
 v_5 + v_6 =  v_7 +v_8 =0 \\
v = \frac12(v_1+ v_3 + v_5 + v_7) \end{array}$ \\\hline
15&$2$& 29&5& 9&5& 5&29&-48& 
$ \begin{array}{c} 2v_1 + 2v_2 + 2v_3 + v_4 + v_5 =0 \\
v =  \frac14(v_1 + 2v_2 + 3 v_3 + 2v_5) \end{array}$ 
                         \\\hline
16&$2$& 29&6& 9&6& 5&29&-48&  $\begin{array}{c}
 v_1 + v_2 + v_3 +v_4 =0 \\
 v_3 + v_4 + v_5 + v_6 =0 \\
v =  \frac14(v_1+ 3v_2 + 2v_4 + 2v_6) \end{array}$
                 
                        \\\hline
\hline\end{tabular}

\end{center}

\section{The Brauer group of toric hypersurfaces}\label{Br}


Let $X$ be a smooth quasi-projective variety over $\C$. The cohomological 
Brauer group $B'(X)$ of $X$ is defined as 
the torsion subgroup in $H^2_{\mathrm{et}}(X, \kO^*_X)$, where 
$\kO^*_X$ denotes the sheaf of units in $\kO_X$. If 
$B(X)$ is the  Brauer group of sheaves of Azumaya algebras on $X$, then 
there exist a canonical injective homomorphism
$B(X) \to B'(X)$. Recently Kresch and Vistoli proved that this 
homomorphism is in fact an isomorphism \cite{KV}. 

Let us summarize some well-known results about the Brauer group.

\begin{theorem} 
The subgroup $B(X)_r \subset B(X)$ of elements $x \in 
B(X)$ with $rx =0$ can be included into the short exact sequence
\[ 0 \to \Pic(X)/r\Pic(X) \to  H^2_{\mathrm{et}}(X,\mu_r) \to B(X)_r \to 0, \]
where $\mu_r$ denotes the subsheaf on $r$-th roots of unity in $\kO_X^*$.
More generally, there is a short exact sequence 
\[ 0 \to 
\Pic(X) \otimes_{\Z}\Q/\Z \to H^2_{\mathrm{et}}(X,\Q/\Z) \to B(X) \to 0,\] 
where  $H^2_{\mathrm{et}}(X,\Q/\Z)$ coincides with the singular 
cohomology group $H^2(X,\Q/\Z)$ in the usual topology of $X$ as a
manifold over $\C$.   
\label{gr1}
\end{theorem}

\prf See \cite{G} II, Theorem 3.1. \hfill $\Box$ 

\begin{coro} If $X$ is a smooth projective variety such that 
$Pic(X) \cong H^2(X, \Z)$, then 
one has the isomorphism 
\[ B(X) \cong {\rm Tors}( H^{3}(X, \Z)). \]
In particular, the last statement holds for smooth 
Calabi-Yau $d$-folds $(d \geq 3)$.
\end{coro} 

\prf Follows from the universal coefficients theorem. \hfill $\Box$

\begin{theorem} 
Let $X$ be a smooth quasi-projective variety and $Y \subset X$ a 
closed subvariety. Then 
the natural homomorphism $ B(X) \to B(X \setminus  Y)$ 
is always injective. This homomorphism  is bijective  if 
$Y$ has codimension $\geq 2$. 
Moreover, if $Y$ is a union of closed irreducible 
subvarieties $Y_1, \ldots, Y_n$ of codimension $1$, then one has the 
exact sequence 
\[ 0 \to B(X) \\\to B(X \setminus  Y) \to \bigoplus_{i =1}^n 
H^1_{\mathrm{et}}(Y_i, \Q/\Z), \] 
 where  $H^1_{\mathrm{et}}(Y_i, \Q/\Z)$ 
 coincides with the singular 
cohomology group $H^1(Y_i,\Q/\Z)$ in the usual topology if $Y_i$ is smooth. 
\label{gr2}
\end{theorem}

\prf See \cite{G} III, Corollary 6.2.  \hfill $\Box$ 
\bigskip

The Brauer group of a $d$-dimensional algebraic torus $\T^d$ over $\C$ 
was computed e.g. by Magid in \cite{Magid}:

\begin{theorem}
Using the canonical isomorphism between the cohomology ring 
$H^*(\T^d, \Z)$ and the exterior algebra $\Lambda^* M$ of the $\Z$-module 
$M$, one has  
the following  canonical isomorphisms
$$B(\T^d) \cong  H^2(\T^d, \Q/\Z) \cong  
\Lambda^2 M \otimes_{\Z} \Q/\Z \cong \Hom( \Lambda^2 N, \Q/\Z). $$
\end{theorem}

The elements in the Brauer group $B(\T^d)$ can be explicitly described
as follows.
Take an integer  $r >1$ and two  characters 
$\alpha, \beta \in M \subset \kO^*_{\T^d}$ 
 of the algebraic torus $\T^d$. 
The {\bf symbol algebra} $(\alpha, \beta)_r$ is generated 
by two elements $x, y$ with relations 
\[ x^r = \alpha,   y^r = \beta, xy = \zeta yx, \]
where $\zeta$ is a primitive $r$-the root of unity. 
If $r =2$, then $(\alpha, \beta)_{2}$ is the classical  quaternionic algebra 
associated with $\alpha, \beta \in \kO^*_{\T^d}$.

Let $B(\T^d)_r$ be $r$-torsion elements in $B(\T^d)$, then 
 $B(\T^d)_r$ is a free $\Z/r\Z$ module with the basis 
\[ (m_i,m_j)_r,  \; \; 1 \leq i < j \leq d, \]
where $m_1, \ldots, m_d$ is a $\Z$-basis of the lattice $M$, i.e. 
\[ B(\T^d)_r = \Lambda^2 (M) \otimes \Z/r\Z.\]
 
\begin{theorem}
Let $\Delta \subset M_{\R}$ be a $d$-dimensional lattice polytope
and let $W \subset \T^d$  be a nondegenerate affine hypersurface defined 
by the equation  
$$f(t) = \sum_{m \in \Delta\cap M} c_m t^m = 0.$$
Then for $d \geq 4$ the natural homomorphism
\[ j^*\; : \; B( \T^d) \to B(W) \]
defined by the embedding $j\, : \, W \hookrightarrow  \T^d$
is an isomorphism.
\label{L-Br}
\end{theorem}

\prf  By \ref{gr1}, we obtain the commutative diagram 
$$\xymatrix{ 0 
\ar[r] & \Pic(\T^d) \otimes_{\Z}\Q/\Z  \ar[d]^{j^*_P} \ar[r] & 
 H^2(\T^d,\Q/\Z)  \ar[d]^{j^*_H} \ar[r] & 
 B(\T^d)  \ar[d]^{j^*}  \ar[r] & 0 \\
0 \ar[r] & \Pic(W) \otimes_{\Z}\Q/\Z   \ar[r] & 
 H^2(W,\Q/\Z) \ar[r] & 
 B(W)   \ar[r] & 0 
}$$
By \ref{L-coh}, $j^*_H$ is an isomorphism. Moreover,  $\Pic( \T^d) =0$ 
since $\T^d$ is a Zariski open subset of $\C^d$. Therefore 
it suffices to prove that 
$\Pic(W)  \otimes_{\Z}\Q =0$ (since by surjectivity of  
$\Pic(W)  \otimes_{\Z}\Q \to \Pic(W)  \otimes_{\Z}\Q/\Z$ would imply 
that   $\Pic(W)  \otimes_{\Z}\Q/\Z=0$). 
Using toric resolution of singularities, one can compactify 
$W$ to a smooth projective 
semiample hypersurface $\widehat{W}$ in some smooth projective        
toric variety $\widehat{\PP}$ such that  $Z:=\widehat{W} \setminus 
W$ is a union of smooth irreducible divisors.
By the vanishing theorem for semiample
line bundles on toric varieties, we obtain that 
\[ h^1(\widehat{W}, \kO_{\widehat{W}}) = \cdots = 
 h^{d-2}(\widehat{W}, \kO_{\widehat{W}}) = 0. \]
In particular, $\Pic(\widehat{W}) \cong H^2(\widehat{W}, \Z)$. 
Therefore, 
$$\dim_{\Q} Pic(W)  \otimes_{\Z}\Q = \dim_{\C} 
 H^2(\widehat{W}, \C) = b_2(\widehat{W}) = b_{2d-4}(\widehat{W}).$$
Let $Z_1, \ldots, Z_k$ be irreducible components of $Z= 
\widehat{W} \setminus 
W$. Using the exact sequence 
\[ \bigoplus_{i =1}^k \Q[Z_i] \to \Pic(\widehat{W}) 
\otimes \Q  \to \Pic(W)  \otimes \Q \to 0, \]
it remains to show that the classes $[Z_1], \ldots, [Z_k] \in  
\Pic(\widehat{W})$ generate $\Pic(\widehat{W}) \otimes \Q$. The latter
follows from the Poincar\'e duality together with  the long exact sequence 
for cohomology with compact supports
\[ \cdots \to H^{2d-4}_c(W, \Q) \to H^{2d-4}_c(\widehat{W}, \Q) \to  
H^{2d-4}_c(Z, \Q) \to \cdots, \]
where $H^{2d-4}_c(Z, \Q) \cong \Q^k$,  $H^{2d-4}_c(\widehat{W}, \Q) 
\cong \Hom ( \Pic(\widehat{W}), \Q)$, because the Hodge component of 
the Hodge type $(d-2,d-2)$ in $H^{2d-4}_c(W, \Q)$ 
is trivial (see \cite{DK}).

\hfill $\Box$  
\bigskip

Let us recall the computation of the  Brauer group of a 
smooth (quasi-projective) toric variety $\PP_{\Sigma}$ obtained  
by Demeyer and Ford \cite{DF}.

By \ref{gr2}, $B(\PP_{\Sigma}) \cong B(\PP_{\Sigma^{(1)}})$ 
where $\Sigma^{(1)}$ is the subfan of all cones of dimension 
$\leq 1$ in $\Sigma$. So we can assume without loss of generality
that $\Sigma = \Sigma^{(1)}$ and $\PP_{\Sigma}$ 
 is a toric compactification of $\T^d$ by divisors 
$D_1^{\circ}, \ldots, D_n^{\circ}$  ($D_i^{\circ} \cong \T^{d-1}$) 
corresponding to lattice vectors $e_1, \ldots, e_n \in N$. 
By \ref{gr2}, we have the exact sequence 
\[ 0 \to B(\PP_{\Sigma}) \to B(\T^d) \to \bigoplus_{i =1}^n 
H^1(D_i^{\circ}, \Q/\Z). \] 
For any $i =1, \ldots, n$ we have the
canonical isomorphisms 
$$H^1(D_i^{\circ}, \Q/\Z) \cong \Hom( N/\Z e_i, \Q/\Z) \cong 
M_{e_i} \otimes \Q/\Z, $$
where 
$$M_{e_i}:= \{ m \in M\; :\; \la x, e_i \ra =0 \}.$$ 
Therefore, one obtains 
\[   B(\T^d) \supset  B(\PP_{\Sigma}) = 
\bigcap_{i =1}^n {\rm Ker}\, \varphi_i,  \]
where $\varphi_i$ denotes the ramification map
\[  \varphi_i\; : \; B(\T^d)  \to H^1(D_i^{\circ}, \Q/\Z) \]
which associates to a given symbol algebra $(\alpha, \beta)_r$ 
its  ramification along the divisor $D_i^{\circ} \subset \PP_{\Sigma}$, 
i.e. a cyclic Galois covering of $D_i^{\circ}$ of degree $r$.   

\begin{theorem} 
Using the canonical isomorphisms  
$$B(\T^d) \cong   \Hom( \Lambda^2 N, \Q/\Z), \; \; 
 H^1(D_i^{\circ}, \Q/\Z) \cong 
\Hom( N/\Z e_i, \Q/\Z),$$ one 
obtains the ramification map $\varphi_i$ from the homomorphism
\[ \psi_i\; : \; N/\Z e_i \to \Lambda^2 N \]
\[ y \mapsto y \wedge e_i, \; \; \forall y \in N \]
by applying the functor $\Hom (*, \Q/\Z)$. In particular, 
the image of $\psi_i$ is the sublattice 
$${\rm Im}\, \psi_i = 
N \wedge e_i \subset  \Lambda^2 N.$$
\label{psi}
\end{theorem}

\prf These statements are contained in   \cite{DF} Lemmas 1.5-1.7. 
\hfill $\Box$ 
\bigskip

Let us define the sublattice  $N': = \sum_{i=1}^n \Z e_i \subset N$ and 
write 
\[ N/N' \cong  \Z/c_1\Z \oplus \cdots \oplus \Z/c_d \Z \]
where $c_1, \ldots, c_d$ non-negative  integers 
such that $c_1 | c_2 | \cdots | c_d$. We can choose a $\Z$-basis 
$y_1, \ldots, y_d$ of $N$ such that $\{ c_iy_i \}_{c_i \neq 0}$ is a 
$\Z$-basis of $N'$. 
Consider  another 
sublattice  
\[  N \wedge N' =  \sum_{i=1}^n 
N \wedge e_i = \sum_{i=1}^n {\rm Im}\, \psi_i \subset \Lambda^2 N. \]
Then the elements $\{ c_i y_i \wedge y_j  \}_{c_i \neq 0, i<j}$ form 
a $\Z$-basis of $N \wedge N'$. Therefore, 
\[  \Lambda^2 N/(N \wedge N') \cong  \bigoplus_{i < j} \Z/c_i \Z. \]
Thus, one obtains (see \cite{DF}, Theorem 1.1):

\begin{theorem}
In the above notations, the Brauer group of a nonsingular toric variety
$\PP_{\Sigma}$ is isomorphic to 
 \[  \bigoplus_{i < j} \Hom( \Z/c_i \Z, \Q/\Z)  . \]
In particular, if $N' \subset N$ is a sublattice of finite index, then 
$\PP_{\Sigma}$ has a finite Brauer group 
\[ B(\PP_{\Sigma})  \cong \Hom( \Lambda^2 N/(N \wedge N'), \Q/\Z) 
\cong  \bigoplus_{i < j} (\Z/c_i \Z)^*. \]
\end{theorem} 

Let $\Delta \subset M_{\R}$ be a $d$-dimensional lattice polytope 
and let $W \subset \T^d$  be a nondegenerate affine hypersurface defined 
by the equation  
$$f(t) = \sum_{m \in \Delta\cap M} c_m t^m = 0.$$
Denote by $\overline{W}$ the closure of $W$ in the projective toric
variety $\PP_{\Delta}$ and by  $\widehat{W}$ a smooth projective
desingularization of  $\overline{W}$ defined by the fan 
$\Sigma \subset N_{\R}$ which is a subdivision of the normal fan 
to $\Delta$. We denote by $N^{(2)}_{\Delta}$ the sublattice in $N$ generated 
by all lattice points $v \in N$ such that the minimum of 
the linear function
$\la *, v \ra$  on $\Delta$ is attained on a face $\theta \subset \Delta$ 
of dimension $\geq 2$. 
Our main result in this section is the following:

\begin{theorem}
In the above notation, for $d\geq 4$ the Brauer group $B(\widehat{W})$ 
of the toric hypersurface $\widehat{W}$ is a cyclic group isomorphic to 
\[ \Hom ( \Lambda^2N/(N \wedge N^{(2)}_{\Delta}), \Q/\Z). \]
\label{main-Br}
\end{theorem}

\prf  
Let $e_1, \ldots, e_n$ be primitive lattice 
generators of $1$-dimensional cones in $\Sigma$ such that 
the minimum of 
the linear function
$\la *, v \ra$  on $\Delta$ is attained on a face $\theta \subset \Delta$ 
of dimension $\geq 1$. Denote by $\Sigma' \subset \Sigma^{(1)}$ the subfan 
of $\Sigma$ consisting of $0$ and all $1$-dimensional cones 
 $\R_{\geq 0} e_1, \ldots,\R_{\geq 0} e_n$. We set  $\widehat{W}^{(1)} = 
\widehat{W} \cap \PP_{\Sigma'}$. Since the complement to $\widehat{W}^{(1)}$ 
in $\widehat{W}$ has codimension $\geq 2$, by \ref{gr2}. it suffices to 
compute the Brauer group $B(\widehat{W}^{(1)})$. By \ref{gr2}, we   
have the commutative diagramm: 
$$\xymatrix{ 0 
\ar[r] & B(\PP_{\Sigma'})  
\ar[d]_{\alpha} \ar[r] &  B(\T^d)   
\ar[d]_{\beta} 
\ar[r]^{\varphi} & \bigoplus_{i=1}^n H^1( D_i^{\circ},\Q/\Z)  \ar[d]_{\gamma} \\
0 \ar[r] & B(\widehat{W}^{(1)}) \ar[r] & B(W) \ar[r]^{\phi} & 
\bigoplus_{i=1}^n  
H^1(Z_i^{\circ}, \Q/\Z)
}$$

By \ref{L-Br}, $\beta$ is an isomorphism. In order to understand the 
homomorphisms $\alpha$ and $\gamma$  
we divide the set of lattice vectors ${\mathcal E}= 
\{ e_1, \ldots, e_n \}$ into 
a union of two disjpoint subsets ${\mathcal E}_1 \cup  {\mathcal E}_2$.
The set  ${\mathcal E}_2$ consists of all $e_i \in {\mathcal E}$ 
such that the minimum of $\la *, e_i \ra$  on $\Delta$ is attained 
on a face $\theta \subset \Delta$ 
of dimension $\geq 2$ and  we define ${\mathcal E}_1 :={\mathcal E} \setminus 
{\mathcal E}_2 $. 

For any $e_i \in {\mathcal E}$ we denote by $\theta_i$
the maximal subface of $\Delta$ on which 
the minimum of $\la *, e_i \ra$  on $\Delta$ is attained.

If  $e_i \in {\mathcal E}_2$, then $d_i:= \dim \theta_i \geq 2$ 
 and the divisor  $Z_i^{\circ}$ is isomorphic to the product   
$V_i \times \T^{d-d_i-1}$, where $V_i$ is a $\theta_i$-
nondegenerate affine hypersurface in a torus $\T^{d_i}$. 
By \ref{L-coh}, 
the ramification homomorphism 
\[ \phi_i \; : \;  B(W) \to H^1(Z_i^{\circ}, \Q/\Z) \]
can be identified with the composition of the 
ramification homomorphism 
\[ \varphi_i \; : \;  B(\T^d) \to H^1( D_i^{\circ},\Q/\Z) \]
with the monomorphism 
\[  H^1( D_i^{\circ},\Q/\Z) \to  H^1(Z_i^{\circ}, \Q/\Z). \]
Thus, we obtain 
\[ B:= \bigcap_{e_i \in {\mathcal E}_2} {\rm Ker}\, \varphi_i = 
 \bigcap_{e_i \in {\mathcal E}_2} {\rm Ker}\, \phi_i. \]
Using  \ref{psi}, we can describe the dual group $B^*$ as cokernel 
of the homomorphisms 
\[ \psi_i \; : \; N/\Z e_i \to \Lambda^2 N, \; \; e_i \in  
{\mathcal E}_2, \]
i.e. we have  
\[ \sum_{e_i \in  {\mathcal E}_2} {\rm Im} \, \psi_i = 
N \wedge  N^{(2)}_{\Delta}\; \; {\rm and } \;\; 
B^* = N/ (N \wedge  N^{(2)}_{\Delta}).   \] 
Therefore, it remains to  prove that  
\[  B =
\bigcap_{e_i \in {\mathcal E}} {\rm Ker}\, \phi_i =  B( \widehat{W}^{(1)}), \]
i.e., that for any   
$e_i \in {\mathcal E}_1$ the image ${\rm Im}\, \phi_i = 
{\rm Im} \, \gamma \circ \psi_i$
is contained in $N \wedge  N^{(2)}_{\Delta}$. Now we consider two cases 
for $e_i \in {\mathcal E}_1$: 

\medskip

\noindent
Case 1: Let $e_i \in  N^{(2)}_{\Delta}$. Then the image 
of $\phi_i$ is contained in $N \wedge e_i \subset N \wedge N^{(2)}_{\Delta}$.

\noindent
Case 2: Let  $e_i \not\in  N^{(2)}_{\Delta}$.  Then the image 
of $\phi_i$ is a  sublattice $N_i \wedge e_i$, where $N_i$ is a direct 
summand of $N$  of rank $d-1$ containing $e_i$ ($N_i$ is a 
sublattice of rank $d-1$ which is orthogonal 
 to the $1$-dimensional face $\theta_i$). The smoothness of $\widehat{W}$
implies that both lattices $N_i$ and $ N^{(2)}_{\Delta}$ contain a 
direct summand of $N$ of rank $d-2$. Therefore, $N_i/  
N_i \cap N^{(2)}_{\Delta}$ is a cyclic group and we have 
$$N_i = \Z v +  N_i \cap N^{(2)}_{\Delta}$$ for some lattice vector $v \in N_i$ 
which proportional to $e_i$. By $v \wedge e_i =0$, we have  
$$ N_i \wedge e_i \subset (\Z v +  N_i \cap N^{(2)}_{\Delta}) \wedge 
e_i \subset  N^{(2)}_{\Delta} \wedge e_i \subset  N \wedge  N^{(2)}_{\Delta}.$$
Thus, all divisors $Z_i$ $(e_i \in {\mathcal E}_1)$ 
do not define non-trivial ramification conditions on the the subgroup $B$ 
of the Brauer group 
$B(W)$ and   
$$ B \cong  B( \widehat{W}^{(1)}) \cong  B( \widehat{W}) \cong
\Hom ( \Lambda^2N/(N \wedge N^{(2)}_{\Delta}), \Q/\Z).$$
\hfill $\Box$

\begin{coro} 
Let $\widehat{W}$ be a smooth projective 
Calabi-Yau $d$-fold  obtained 
as a crepant desingularization of a nondegenerate 
projective hypersurface $\overline{W} \subset \PP_{\Delta}$ 
corresponding to a reflexive polytope 
$\Delta$ of dimension $d \geq 4$.
Then the Brauer group   $B(\widehat{W})$  is a cyclic group isomorphic to
\[ \Hom ( \Lambda^2N/(N \wedge N_{\Delta^*}''), \Q/\Z), \]
where $N_{\Delta^*}'' \subset N$ is the sublattice generated by 
all lattice points in $N \cap \Delta^*$ which are contained in 
 faces of $\Delta^*$ of codimension $> 2$.
\end{coro} 

\prf The statement follows immediately form \ref{main-Br} and 
from the easy combinatorial equality $N_{\Delta}^{(2)} = N_{\Delta^*}''$ 
using the duality between faces of dual reflexive polytopes 
$\Delta \subset M_{\R}$ and $\Delta^* \subset N_{\R}$. \hfill $\Box$

\section{The Brauer group of Calabi-Yau 3-folds}\label{Br2}


\bigskip

Let $\Delta$ be a $4$-dimensional reflexive polytope. Then 
\[ N_{\Delta^*}'':= \sum_{v \in \theta^*\cap N, \dim \theta^* = 1} \Z v 
\subset N\]
is a sublattice of finite index. We write 

\[ N/N_{\Delta^*}'' = 
\Z/c_1\Z \oplus  \Z/c_2\Z \oplus \Z/c_3\Z \oplus  \Z/c_4\Z, \]
where $c_1 | c_2 | c_3 | c_4$. In this case  $c_1 = c_2 = 1$ and 
\[ B(\widehat{W}) \cong \Hom( 
\Lambda^2 N/(N \wedge N_{\Delta^*}''), \Q/Z)  \cong  
\Z/c_3\Z. \]

Using the Calabi-Yau date \cite{KS2},  one can check all 
 473 800 776 reflexive
polytopes $\Delta$ and find exactly 16 reflexive polytopes of dimension 4 
with $c_3 > 1$. For each of these 16 polytopes, one obtains 
\[  N/N_{\Delta^*}'' \cong   \Z/p\Z\oplus \Z/p\Z, \;\; p = 2,3,5. \]
So we have again 16 families with the  cyclic Brauer group $\Z/p\Z$. 
The main observation is that these 16 families are exactly the mirrors
of 16 non-simply connected families of Calabi-Yau 3-folds 
with the cyclic fundamental group $\Z/p\Z$ $(c =2,3,5)$ (see Section \ref{ex}).
Therefore, we can exchange the lattices $N \leftrightarrow M$ and 
polytopes  $\Delta^* \leftrightarrow \Delta$ so that for all $4$-dimensional 
reflexive polytopes we have  the 
isomorphisms
\begin{equation}
\Lambda^2 M/(M \wedge M_{\Delta}'') \cong 
 N/N_{\Delta^*}' \;\; {\rm and} \;\;  
\Lambda^2 N/(N \wedge N_{\Delta^*}'') \cong 
 M/M_{\Delta}' 
\label{comb-dual}
\end{equation}
Unfortunately we do not have any natural 
combinatorial explanation of the
isomorphisms (\ref{comb-dual}). 
\bigskip

The sublattice $N_{\Delta^*}'' \subset N$ defines a (ramified) Galois 
covering $\widehat{W}'' \to \widehat{W}$ with the Galois group 
$G = N/N_{\Delta^*}'' \cong  \Z/p\Z\oplus \Z/p\Z$. 
If we consider $G$ as a $2$-dimensional vector space over the 
finite field $\F_p$, then 
the Brauer group $B(\widehat{W})$ is isomorphic to 
$\Hom (\Lambda^2 G, \Q/\Z)$, i.e, to a $1$-dimensional $\F_p$-space.  
This allows to look at  $B(\widehat{W})$ as the discrete torsion 
group $H^2(G, U(1))$  for the orbifold  $\widehat{W}''/G$.  
\bigskip 

We finish this section by describing 
a typical example $\widehat{W}^*_1$  of a Calabi-Yau $3$-fold from our list 
(the mirror of $\widehat{W}_1$ from Section \ref{ex}) 
with Brauer group $B(\widehat{W}^*_1) \cong   \Z/5\Z$. 
This example gives a  $1$-parameter family of Calabi-Yau 3-folds 
$\widehat{W}^*_1$ 
having a close connection to the counterexample to the global 
Torelli considered by Aspinwall-Morrison and B. Szendr\H{o}i  
\cite{AM1,Sz}.  

Consider the action of the group 
$G:= \mu_5 \times \mu_5$ on $\PP^4$ by 
\[ \zeta_1 \to (1, \zeta_1,\zeta_1^2,\zeta_1^3,\zeta_1^4) \]
\[ \zeta_2 \to (1, \zeta_2,\zeta_2^3,\zeta_2,1). \] 
This action is not free and the quotient
$\PP_{\Delta_1^*} := \PP^4 /G$ has cyclic singularities of type 
$\frac{1}{5}(3,1,1)$ along 10 curves.
Let $\widehat{W}^*_1$ be the crepant desingularization of a generic 
Calabi-Yau hypersurface $\overline{W}_1^* \subset \PP_{\Delta^*_1}$. Then   
\[ \pi_1(\widehat{W}^*_1) = 0 =  B(\widehat{W}_1) \]
\[ B(\widehat{W}^*_1) = \Lambda^2 G \cong  \mu_5 \cong  
\pi_1(\widehat{W}_1). \]

\bigskip


\end{document}